\documentclass[twoside,12pt]{amsart}

\usepackage[english]{babel}
\usepackage{amsfonts}
\usepackage{amscd}
\usepackage{amsmath}

\newtheorem{Prop}{Proposition}
\newtheorem{Cor}[Prop]{Corollary}

\newtheorem{Def}[Prop]{Definition}

\newtheorem{rmk}[Prop]{Remark}
\newenvironment{Rmk}{\begin{rmk}\em}{\end{rmk}}
\newtheorem{exm}[Prop]{Example}

\newtheorem{prf}{Proof}

\newenvironment{Prf}{\begin{prf}\em}{\qed\end{prf}}
\newtheorem{prff}{}

\newcommand{\NOT}[1]{}
\newcommand{\pa}{\par\medskip}

\newcommand{\EM}{\em\bf}
\newcommand{\BAR}[1]{\overline{#1}}
\newcommand{\HAT}[1]{{\buildrel{\,\,\wedge}\over{#1}}\!}

\newcommand{\eps}{\varepsilon}

\newcommand{\al}{\alpha}
\newcommand{\be}{\beta}
 \newcommand{\ga}{\gamma}

\newcommand{\lr}{{lim-rim}}
\newcommand{\Lr}{{Lim-Rim}}
\newcommand{\U}{{^\star}}
\newcommand{\Uz}{{{}^\star_0}}
\newcommand{\T}{{^\star}\!\!}
\newcommand{\NS}{non-standard}
\newcommand{\NNS}{Non-Standard}
\newcommand{\NA}{non-standard analysis}
\newcommand{\NNA}{Non-Standard Analysis}
\newcommand{\Dem}{\EM}
\newcommand{\True}{\textbf{True}}
\newcommand{\False}{\textbf{False}}
\newcommand{\TF}{\{\True,\False\}}
\newcommand{\Hom}{\text{Hom}}
\newcommand{\cd}[1]{\mathfrak{#1}}

\newcommand{\bbb}[1]{\mathbb{#1}}

\newcommand{\bN}{\bbb{N}}
\newcommand{\bQ}{\bbb{Q}}
\newcommand{\bR}{\bbb{R}}

\newcommand{\cE}{\mathcal{E}}
\newcommand{\cF}{\mathcal{F}}
\newcommand{\cH}{\mathcal{H}}
\newcommand{\cI}{\mathcal{I}}
\newcommand{\cL}{\mathcal{L}}
\newcommand{\cN}{\mathcal{N}}
\newcommand{\cP}{\mathcal{P}}
\newcommand{\cS}{\mathcal{S}}
\newcommand{\cT}{\mathcal{T}}
\newcommand{\cU}{\mathcal{U}}
\newcommand{\cV}{\mathcal{V}}

\newcommand{\AP}{as per Definition \ref{Def:na}}
\newcommand{\BC}{$|B|$-confined \NA\ induced by}

\newcommand{\cSET}{\cS\cE\cT}
\newcommand{\cFIN}{\cF\cI\cN}

\title{Non Standard Analysis as a Functor, as Local, as Iterated}
\author{Eliahu Levy}
\address{Department of Mathematics,
Technion -- Israel Institute of Technology,
Haifa 32000, Israel}
\email{eliahu@techunix.technion.ac.il}

\date{}

\begin{document}

\begin{abstract}
This note has several aims. Firstly, it portrays a \NA\ as a functor, namely a functor $\U$ that maps any set $A$ to the set $\U A$ of its \NS\ elements. That functor, from the category of sets to itself, is postulated to be an equivalence on the full subcategory of finite sets onto itself and to preserve finite projective limits (equivalently, to preserve finite products and equalizers). Secondly, ``Local'' \NA\ is introduced as a structure which we call \lr, in particular exact \lr s. The interplay between these, and ultrafilters and ultrapowers, and also cardinality relations and notions depending on a cardinality such as saturation and what we call ``confinement'' and ``exactness'', are investigated.

In particular, one constructs \NS\ analyses, with ``good'' kinds of \lr. In these one may say that $\U A$ -- ``the adjunction of all possible limits from $A$'' -- plays a role analogous to that of the algebraic closure of a field -- ``the adjunction of all roots of polynomials''. Then in the same spirit as with the latter, one has uniqueness up to isomorphism, and also universality and homogeneity, provided one has enough General Continuum Hypothesis. The cardinality of $\U A$ will be something like $2^{2^{|A|}}$ -- the same as that of the set of ultrafilters in $A$, and one has a high degree of saturation.

Also, one notes that the functor $\U$ can be applied again to $\U A$, giving $\U\U A$,
$\U\U\U A$, and so forth. In particular, we focus on the two different embeddings of $\U A$ into $\U\U A$ and prove some of their properties, with some applications.
\end{abstract}

\maketitle

\tableofcontents

\setcounter{section}{-1}

This note is a sequel to the less formal \cite{Levy}. Here I make things more precise, clarify and develop them, and add some material about iterated \NA. In particular, a \NA\ (which here gives ``non-standard elements'' to any set) is approached as a functor (a connection between \NA\ and category theory was given in \cite{Saigo}). ``Local'' \NA\ is introduced as a structure which I call \lr, in particular exact \lr s. The interplay between these, and ultrafilters and ultrapowers, and also cardinality relations and notions depending on a cardinality such as saturation and what I call ``confinement'' and ``exactness'', are investigated.\pa
 
In the second section, dedicated to iterated \NA, the two embeddings $\nu$ and $\U\nu$ of $\U A$ into $\U\U A$ are focused upon. These embeddings are highly distinct. In the last paragraph a simple criterion (Prop.\ \ref{Prop:nunustar}) that the values they take satisfy a given standard relation is given, which has some interesting consequences and applications.

\section{Notations}
Denote the cardinality of a set $A$ by $|A|$.
For a cardinality $\cd{m}$, denote by $\cd{m}^+$ the successor cardinality.
Denote the power set of a set $A$ (the set of all subsets of $A$) by $\cP A$.

\section{\NNA}

\subsection{\NNA\ as a Functor}
\begin{Def}\label{Def:na}
Let $\cSET$ be the category of sets and maps between them, and let $\cFIN$ be its
full subcategory of finite sets and maps between them.
A {\Dem\NA} will be defined as a (covariant) functor $\U:\cSET\to\cSET$ which (i) maps $\cFIN$ to itself as a self-equivalence of $\cFIN$; (ii) preserves finite projective limits.
\end{Def}

$\U$ of a set $A$ will be denoted by $\U A$ and its elements will be referred to as the {\Dem\NS\ elements} of $A$ (thus, if $A$ is, say, the set of integers, then the elements of $\U A$ will be called $\U$integers). Similarly the image by $\U$ of a map $f$ will be denoted by $\U f$.\pa

In particular, for every $a$ the set $\{a\}$ is a singleton, thus an object of $\cFIN$, and $\U$ being an equivalence on $\cFIN$, $\U\{a\}$ is also a singleton and its sole element will be denoted by $\nu(a)$ ($\nu$ from ``new''). In many cases, we will tend to identify $a$ with $\nu(a)$ (see below).\pa

As for condition (ii), it is equivalent to $\U$ preserving finite products (equivalently products of two sets $A$ and $B$, being the projective limit of the diagram $A,B$ (no maps)\NOT{, characterized by the appropriate universal property}) and equalizers $\{c\in C|f_1(c)=f_2(c)\}$ of two maps $f_1,f_2:C\to D$, (i.e.\ the projective limit of the diagram $C{\to\atop\to}D$ -- the maps are $f_1$ and $f_2$\NOT{ -- characterized by the appropriate universal property}).\pa

This means that the images by $\U$ of the projections $A\times B\to A,\,A\times B\to B$ make $\U(A\times B)$ the product $\U A\times\U B$; and if $E$ is the subset $\{c\in C\,|\,f_1(c)=f_2(c)\}$, i.e., with the inclusion map $E\to C$, the equalizer of $f_1$ and $f_2$, then $\U$ of the inclusion map makes $\U E$ the equalizer of $\U f_1,\U f_2:\U C\to\U D$. But since any subset of $C$ is an equalizer, this means that $\U$ of any inclusion map $E\to C$ is injective, and can be viewed as an inclusion $\U E\subset\U C$, and in the above situation $\U E=\{\U c\in\U C\,|\,\U f_1(\U c)=\U f_2(\U c)\}$.\pa

In particular, $\U$ of the diagonal of $A\times A$ -- the equalizer of the two projections -- can serve as the diagonal of $\U A\times\U A$.\pa

\begin{exm}\label{exm:na}

Of course, a trivial example for \NA\ is the identity functor, where ${\U}A=A$ and all \NS\ elements are, in fact, standard -- of the form $\nu(a)$.\pa

To get a non-trivial example, Let $S$ be a fixed set of indices and $\cU$ a fixed ultrafilter in $S$. For any set $A$ we have the {\Dem ultrapower} $A_\cU:=A^S/\cU$, defined as the power $A^S$ modulo the equivalence relation: equality modulo $\cU$.\pa

We have here a covariant functor $A\mapsto A_\cU$ (which depends on $S$ and $\cU$)
from the category of sets to itself.  The image $f_\cU$ by the functor of a map $f:A\to B$ between sets will map $[\xi]\in A_\cU$, where $\xi:S\to A$, to $[f\circ\xi]\in B_\cU$.\pa

This ultrapower functor satisfies (i) and (ii) of Definition \ref{Def:na}. Indeed, (i) is immediate, because if $F$ is finite then for any $\xi:S\to F$ there is a unique $a\in F$ such that $\xi^{-1}(\{a\})\in\cU$, and then $\xi$ is, modulo $\cU$, the constant $a$. As for (ii), this functor preserves finite products - $(A\times B)_\cU$ is identified with
$A_\cU\times B_\cU$ in an obvious manner, and if $f_1,f_2:C\to D$, and
$E=\{c\in C\,|\,f_1(c)=f_2(c)\}$ is the equalizer, then the equalizer of $(f_1)_\cU,(f_2)_\cU:C_\cU\to D_\cU$ is
\begin{eqnarray*}
&&\{[\xi]\in C_\cU\,|\,(f_1)_\cU([\xi])=(f_2)_\cU([\xi])\}=\\
&&=\{[\xi]\in C_\cU\,|\,[f_1\circ\xi]=[f_2\circ\xi]\}=\\
&&=\{[\xi]\in C_\cU\,|\,\exists U\in\cU,\,\,(f_1\circ\xi)|_U=(f_2\circ\xi)|_U\}=\\
&&=\{[\xi]\in C_\cU\,|\,\exists U\in\cU\,\,\xi(U)\subset E\}=\\
&&=\{[\xi]\in C_\cU\,|\,\exists\eta:S\to E\,\,[\xi]=[\eta]\},
\end{eqnarray*}
which is identified with $E_\cU$.\qed\pa

This means that taking $\U A:=$ the ultrapower $A_\cU$, gives us a \NA, \AP.\pa
\end{exm}

\subsection{\NNA\ as a Functor, continued}\label{S:NAFco}
If $A$ is a set and $a\in A$, there is an inclusion map $\{a\}\to A$, and its $\U$ is an injection $\U\{a\}\to\U A$, i.e.\ $\{\nu(a)\}\to\U A$, identifying $\nu(a)$ with an element of $\U A$ (Consequently, if $A\ne\emptyset$ then $\U A\ne\emptyset$) -- if $a$ belongs to some (standard) set -- has some property -- then $\nu(a)$ is one of the \NS\ elements of the set -- has the $\U$property. These $\nu(a)$'s, with $a\in A$, might be referred to as the {\Dem standard elements} of $A$, but (if $A$ is infinite) $A$ may have other elements, not of the form $\nu(a)$,\,\,$a\in A$.\pa

Similarly, $\nu(a_1,\ldots,a_n)$ is identified with $(\nu(a_1),\ldots,\nu(a_n))$, and if $a_1\,\ldots,a_n$ satisfy some relation $R\subset A_1\times\cdots\times A_n$ then $\nu(a_1),\ldots,\nu(a_n)$ satisfies $\U R$.\pa

We will sometimes loosely view any finite set $F$ as a set of a kind of {\em truth-values} and any map between finite sets $f:F_1\to F_2$ as a kind of {\em logical operation}, \`a la {\bf AND}, {\bf OR}, {\bf NOT}. Then the equivalence that $\U$ induces on $\cFIN$ is used as an identification for these truth-values and logical operations. Maps $f:A_1\times A_2\times\cdot\times A_n\to F$, where $A_i, i=1\ldots n$ are sets and $F$ is finite, are viewed as {\em predicates} on $A_1,A_2,\ldots,A_n$. Composing them with maps between finite sets defines a kind of {\em logical operations among predicates}. But the fact that $\U$ preserves finite products and equalizers makes, e.g., a $\TF^2$-valued predicate the same as a pair of $\TF$-valued predicates, and if $F'\subset F$ then $F'$-valued predicates embed into $F$-valued ones. Therefore the $\TF$-valued predicates determine all others, while applying $\U$ to predicates commutes with all ``logical operations'', in other words, applying $\U$ to $\TF$-valued predicates commutes with the propositional calculus operations {\bf AND}, {\bf OR}, {\bf NOT} etc.\pa

As said above, we identify $\nu(\True)=\True$,\,\,$\nu(\False)=\False$.\pa

For any $\TF$-valued predicate $f$, if $E$ is the subset of
$A_1\times A_2\times\cdots\times A_n$ where it holds, then $E$ is the equalizer of $f$ with the constant map $A_1\times A_2\times\cdots\times A_n\to\True$, the latter being the composition of $\{\True\}\to\TF$ with the sole map
$A_1\times A_2\times\cdots\times A_n\to\{\True\}$. Hence its $\U$ is the constant
${\U}A_1\times{}{\U}A_2\times\cdots\times{}{\U}A_n\to\True$ and $\U E$, as a subset of
${\U}A_1\times{}{\U}A_2\times\cdots\times{}{\U}A_n$, will be its equalizer with $\U f$, i.e.\ the set where $\U f$ holds.\pa

As for quantifiers, we shall show that the operation $\U$ on $\TF$-valued predicates commutes with them too. Since it commutes with {\bf NOT}, it suffices to show that it commutes with $\exists$.\pa

Given a predicate $f:A_1\times A_2\times\cdots\times A_n\to\TF$, to which we apply $(\exists x_1\in A_1)$, we get the predicate $(\exists x_1\in A_1)f$, and we have to prove that the $\TF$-valued functions (``predicates''):
\begin{equation}\label{eq:exists}
\U[(\exists\,x_1\in A_1)\,f(x_1,\ldots,x_n)]=
(\exists\,\U x_1\in{}{\U}A_1)\,\U f(\U x_1,\ldots\U x_n).
\end{equation}
If $A_1=\emptyset$ this is immediate. Thus suppose $A_1\ne\emptyset$, consequently $\U A_1\ne\emptyset$.\pa

Firstly, denote $g(x_2,\ldots,x_n):=(\exists x_1\in A_1)f(x_1,\ldots,x_n)$.
Then $A_1\times\{g=\True\}\supset\{f=\True\}$. Hence ${\U}A_1\times\{\U g=\True\}\supset\{\U f=\True\}$, implying that ${\U}A_1\times\{\U g=\True\}$,
i.e.\ the product of ${\U}A_1$ with the LHS of (\ref{eq:exists}),
contains the product of ${\U}A_1$ with the RHS of (\ref{eq:exists}).
Since ${\U}A_1\ne\emptyset$, this implies LHS $\supset$ RHS.\pa

For the reverse inclusion we make use of a map $j$, defines on the subset $g=\True$ and mapping each such $(x_2,\ldots,x_n)$ to some $x_1\in A_1$ such that
$f(x_1,x_2,\ldots,x_n)$ holds. (Note that if $n=1$,\,\,$A_2\times\cdots\times A_n$ is a singleton -- a product of an empty family.) $\U\{g=\True\}$ is the LHS of (\ref{eq:exists}) and taking for any $(\U x_2,\ldots,\U x_n)\in\,\U\{g=\True\}$,\,\,$\U x_1:=\U j(\U x_2,\ldots,\U x_n)$, proves that $\U\{g=\True\}$ is contained in the RHS.\qed\pa

Therefore any \NA\ \AP\ satisfies the {\em Transfer Principle}: predicates (in the usual sense, i.e.\ $\TF$-valued), commute with all the logical operations of a first-order theory, where quantifiers taken over a set $A$ should correspond to quantifiers taken over $\U A$. In particular, a sentence which is a first-order logical combination of standard sets $A$ is true if and only if the sentence which is the same first-order combination, each $A$ replaced by ${\U}A$, (hence each element $a$ by $\nu(a)$), is true.\pa

Since, by the above, inclusion of sets implies inclusion of the $\U$-sets,%
\footnote{One may, say, argue thus: if $C_1$,\,$C_2$ are subsets of a set $A$, with inclusion maps $i_1:C_1\to A$,\,$i_2:C_2\to A$, then $C_1\subset C_2$ if and only if there is an injection $j:C_1\to C_2$ such that $i_1=i_2\circ j$.}
we may define the $\U$ of a class (e.g.\ of the relation $\in$ between sets) as the direct limit of the $\U$ of all sets contained in the class.\pa

Note that often we shall use the same notation for a relation $R$ and for $\U R$, for example write $<$ for $\T<$ or $\in$ for $\T\in$. By Transfer, any first-order property of $R$ holds for $\U R$. Thus, for example, if $\U S$ is a $\U$set (i.e.\ a \NS\ element of some family of sets), then $\U S$ is determined by the set of $\U x$ such that $\U x\,\,\T\in\U S$, the latter set called the {\Dem scope} \cite{MachoverHirschfeld} of $\U S$ and sometimes denoted by $\HAT{\U S}$\,. Thus $\U x\,\,\in{}\HAT{\U S}\Leftrightarrow\,\U x\,\,\T\in\,\U S$. But sometimes we shall identify $\U S$ with $\HAT{\U S}$\,, view $\T\in$ as a particular case of $\in$, and thus view a $\U$set as a set. Yet, as is well known, not all sets of \NS\ elements are of the form $\HAT{\U S}$\,\, for some $\U$set $\U S$. Sets of \NS\ elements of the latter form are called {\Dem internal} and by our above identifications we may view internal sets as the same as $\U$sets. Sets of \NS\ elements which are not internal will be called {\Dem external}.%
\footnote{\label{f:i}For example, by Transfer any bounded internal set of positive $\U$integers has a greatest element, and any positive $\U$integer less or equal than some standard $n$ (i.e.\ $\nu(n)$)\,\,$\U$belongs to $\U\{1,\ldots,n\}=\{1,\ldots,n\}$, i.e.\ is a standard $\nu(m),\, m\le n$. Therefore if $\U\bN$ contains any element which is not standard, then this element is bigger than all the standard positive integers, the set of these standard positive integers is bounded but has no greatest element, hence cannot be internal, i.e.\ is external.}

Note, that for a set $A$,\,\,
$$\forall\,x\,\,(x\in A\Leftrightarrow\exists\,y\in\{A\}\,\,x\in y)$$
is true, therefore
$$\forall\,\U x\,\,(\U x\in\U A\Leftrightarrow
\exists\,\U y\in\U\{A\}\,\,\U x\,\,{\T}\in\U y).$$
is true, where we took the $\U$ of the relation $\in$ between elements and sets. But $\U\{A\}=\{\nu(A)\}$, hence we find that as sets $\U A=\nu(A)$ (i.e. $\U A=\HAT{\nu(A)}$). But we shall keep using different notations: $\U A$ where we emphasize the set of (\NS) elements, and $\nu(A)$ when we emphasize the single element.\pa

\subsection{Features of the Conventional Approach to \NNA}
\NA\ was constructed by its founder, Abraham Robinson
(\cite{Robinson}) as a \NS\ model of part of Set Theory. In his and
following treatments (see several treatises in the bibliography)%
\footnote{In another direction, beginning with E.~Nelson \cite{Nelson}, one
constructs a special Set Theory for \NA.}
the idea is to start with a ``big'' structure -- a set $V$, with the restriction to it of the membership relation $\in$, which will be big enough so that usual mathematics could be done in it (should contain, e.g.\ $\bN$ (the natural numbers), $\bR$ (the real numbers), sets of them, sets of these etc.) $V$ is to be provided with a ``mirror'' ``\NS''
structure $\U V$ with a relation $\T\in$, viewed as an ``interpretation'' of the $\in$ of $V$ (Also, ``interpret'' ``$=$'' in $V$ by ``$=$'' in $\U V$), so that

\begin{enumerate}
\item
$V$ is embedded in $\U V$ -- identified with a subset of $\U V$. The members of
$V$ (thus considered as members of $\U V$) are referred to as {\em standard}
(yet by ``\NS'' elements we mean any members of $\U V$).

\item
A {\em Transfer Principle} holds: any (in some constructions one has to say:
``bounded'') first-order logical expression with members of $V$ as constants which is a true sentence when quantifiers are interpreted relative to $V$ will remain a true sentence if (the constants -- members of $V$ -- are replaced by their identified images in $\U V$ and) $\in$ is replaced by $\T\in$ and the quantifiers are interpreted relative to $\U V$. In this sense $\U V$ is a model to the true in $V$ sentences of the first-order theory with $\in$ and all elements of $V$ as constants.

\item
By Transfer, for any standard set $A$, a standard $a\in V$ will be an $\T\in$-member of it if and only if it is a usual member, but in $\U V$\,\,$A$ may have other, {\em\NS} members. The set of all elements of $\U V$ which $\T\in$-belong to $A$ is denoted by $\U A$. Again by Transfer, $\U (A\times B)$ is canonically identified with $\U A\times\,\U B$, thus for relations, functions etc.\ we have the $\U$-relation, $\U$-function etc.\ among members of $\U V$ (which for standard elements coincides with the original relation, function
etc.) When there is no danger of misunderstanding, one often uses the same
symbol (e.g.\ $+$, $<$) for the $\U$-relation etc.\ as for the original.

\item
$\U V$ is indeed ``\NS'' -- it contains elements not in $V$, and one has features, usually some kind of ``saturation'' (see below), which will ensure, say, that there is an $\U n\in\U V$ which $\T\in$-belongs to the set $\bN$ of natural numbers, but such that for every standard $m\in N$,\,\,$m\;\T<\U n$ (i.e.\ $(m,\U n)$ $\T\in$-belongs to the graph of $<$). Such $\U n$'s are naturally referred to as
``infinite'' or ``unlimited''%
\footnote{In fact, as in footnote \ref{f:i}, for $\bN$ the existence of ``infinite'' members needs only the existence of some not standard members -- one easily proves using Transfer that all $\U$-members of a finite standard $A$ are standard, hence any $\U$-member of $\bN$ which is not standard is automatically bigger than all standard
numbers.}
and by Transfer they have reciprocals in $\U\bR$, referred to as ``infinitesimal''.
\end{enumerate}

The ``philosophy'' here is to develop the theory from these properties, paying
less attention to the way such $\U V$ is constructed (or proved to exist).
The usual construction is by an ultrapower or a modification of it.
Anyway, one has to have a set $V$ to make the construction, hence only
subsets of it will have \NS\ members.
\pa

\subsection{Ultrafilters of \NNS\ Elements}
Let us be given a \NA, \AP. For any set $A$, and any
$\U x\in\U A$, and for any $S\subset A$, we have $\U S\subset\U A$, thus either
$\U x\in\U S$ or $\U x$ is not in $\U S$ $\Leftrightarrow$ it is in the complement
${\U}A\setminus\U S=\U (A\setminus S)$. Put otherwise: for any $\TF$-valued predicate
$f:A\to\TF$ we have $\U f:{}\U A\to\U\TF=\TF$ and $\U f(\U x)$ is either $\True$ or $\False$.\pa

Hence any $\U x\in{}\U A$ induces a map $\TF^A\to\TF$, which is easily seen to be a Boolean-algebra homomorphism, i.e.\ an ultrafilter in $A$.\pa

In fact, we may replace $\TF$ by any other finite set $F$, and view the set $\be A$ of the ultrafilters in $A$ as the set of natural mappings between the functor from $\cFIN$ to $\cSET$\,\,$F\mapsto F^A=\Hom(A,F)$, and the inclusion functor $\cFIN\to\cSET$.\pa

Then to a $\U x\in\U A$ corresponds the ultrafilter: the mapping $[f\mapsto\U f(\U x)]:\Hom(A,F)\to F$. This is clearly a natural mapping between the functors $A\mapsto\U A$ and $A\mapsto\be A$.\pa

Note, that although the functor $A\mapsto\be A$ is a canonically defined functor on sets, which maps injections to injections, it does not preserve finite products (hence is not a \NA\ \AP). The projections $A\times B\to A$ and $A\times B\to B$ define onto mappings $\be(A\times B)\to\be A$ and $\be(A\times B)\to\be B$, hence induce a $\be(A\times B)\to\be A\times\be B$, which is also onto, but is usually not one-to-one. Thus this functor does not satisfy (ii) of Definition \ref{Def:na} (while it satisfies (i)).\pa

One may view a \NA\ as a quest to extend sets by adjoining all (or some) possible limits. This in close analogy in spirit to extending a field to its algebraic closure by adjoining all possible roots of algebraic equations. In the latter case one would like to simply define the adjoined elements as labeled by the irreducible polynomial they satisfy, but that is hampered by such facts as the polynomial whose roots are the sums of the roots of two irreducible polynomials being not irreducible.%
\footnote{Which may be viewed as the ``reason'' for the existence of several
conjugate roots: the polynomial $d(t)$ whose roots are the differences of the
roots of the irreducible $f(t)$ has a factor $t$ but also other irreducible
factors, testifying to differences $\ne0$.}
The algebraic closure thus contains conjugate elements with the same irreducible polynomial and for any such elements there is an automorphism of the algebraic closure exchanging them.\pa

Similarly in our case we would like to label the adjoined ``\NS\ elements'' by
the ultrafilters on the original set. After all, these ultrafilters (and ultrafilters in finite products of sets) give all the possible ways properties (and relations) can conceivably hold or not. But that would not work precisely because the functor $\be$ does not preserve finite products. There is a naturally defined Cartesian product of filters, but the Cartesian product of two ultrafilters $\cU\in\be A,\cV\in\be B$ is in general not an ultrafilter -- it is the intersection of all ultrafilters in $A\times B$ that map to $(\cU,\cV)$ in the above mapping $\be(A\times B)\to\be A\times\be B$.\pa

We too will have a construction (\S\ref{S:car}) where there will be ``conjugate'' \NS\ elements with the same ultrafilter and any such conjugate elements will be exchangeable by an automorphism, provided one has enough General Continuum Hypothesis (GCH), and then, in the same spirit as with algebraic closures, one has uniqueness (for fixed ``basis'' $B$ -- see below) up to isomorphism, and also universality and homogeneity.
\pa

\subsection{\Lr s -- Local \NNA}
Given a \NA, one may focus on a particular set $B$ and try to extract what the \NA\ says about this $B$ (hence about its subsets).\pa

The information telling which \NS\ elements satisfy which relations is encoded in the ultrafilters on $B^I$,\,$I$ finite. The \NA\ gives us a natural mapping $(\U B)^I\to\be(B^I)$ between these two (contravariant) functors on $I$ from $\cFIN$ to $\cSET$. If for any set $I\in\cSET$ we will define $\be_c(A^I)$ as the set of ultrafilters on the {\em cylinder Boolean algebra of $B^I$} (the members of the cylinder Boolean algebra are the subsets of $B^I$ that depend only on a finite number of coordinates. In case $I$ is finite this is just the power set of $B^I$), then any natural $(\U B)^I\to\be(B^I)$ on the finite sets $I$ extends in a unique way to a natural $(\U B)^I\to\be_c(B^I)$ on all sets $I$. Moreover, we can take $I=\U B$ itself and then the image by the natural mapping of any element of $(\U B)^I$ will be known if we know the image in $\be_c\left( B^{{}{\U}B}\right)$ of the identity family $\in(\U B)^{\U B}$.

\begin{Def}
We say that a set $E$ has the structure of a {\Dem\lr} over a set%
\footnote{we shall usually have $B$ infinite.}
(basis) $B$ if either of the equivalent following is given. (i) An  ultrafilter $\cL$ on the cylinder Boolean algebra of $B^E$, which will be referred to as the {\Dem defining cylinder ultrafilter} of the \lr\ (then we shall also speak of the \lr\ $(E,\cL)$ over $B$). (ii) A natural mapping between the (contravariant) functors on sets $I$: $E^I$ and $\be_c(B^I)$, where the latter denotes the set of ultrafilters in the cylinder Boolean algebra. It will be referred to as the {\Dem defining natural mapping} of the \lr. (Note that these functors and the natural mapping are defined for any $I$, yet they are determined by giving them for the finite $I$.)
\end{Def}

Here, the defining cylinder ultrafilter $\cL$ in $B^E$ can be recovered as the image, by the natural mapping, in $\be_c\left( B^E\right)$, of the identity family in $E^E$. On the other hand, if $\cL$ is given, what it does is deciding, for any (finite) family $\eta\in E^I$ indexed by a finite $I$ (pushing $\cL$, using the map $B^E\to B^I$ induced by $\eta:I\to E$, to
an ultrafilter in $B^I$) whether relations -- subsets of $B^I$ -- hold or
not, thus transferring such $I$-relations in $B$ (= subsets of $B^I$) into
$I$-relations in $E$, as a \NS\ setting should. This new $I$-relation in $E$ will be referred to as the {\em transfer} or $\U$ of the original $I$-relation in $B$.\pa

To put it otherwise: for an $n$-tuple $(e_1,\ldots,e_n)$ of elements of $E$, a
relation $R\subset B^n$\,\,$\U$-holds for $(e_1,\ldots,e_n)$ if and only if for $\psi\in B^E$,
$(\psi_{e_1},\ldots,\psi_{e_n})\in R$ holds modulo $\cL$.\pa

To conclude:\pa

\begin{Prop}
Given a \NA\ \AP, it induces on any set $B$ a \lr\ $E$ over $B$, where $E=\U B$.
\end{Prop}
And as said above, this \lr\ captures which relations about $B$, or its subsets, are satisfied by $n$-tuples of \NS\ elements $e_1,\ldots,e_n$. We may say that a \lr\ $E$ over $B$ defines a ``local'' \NA, with $E$ the set of \NS\ elements of $B$, completely satisfactory as a \NA\ as long as one restricts oneself only to $B$ or its subsets.\pa

\begin{Rmk}
There can be sub-\lr s of a \lr\ $E$ over $B$ in two ways (or both combined):
Firstly, any subset $E'\subset E$ has the structure of a \lr\ over $B$,
and secondly for any $B'\subset B$ the set $\U B'$ of the $e\in E$ that satisfy
the $\U $ of `` belongs to $B'$ '' forms a \lr\ over $B'$. Instead of
talking about sub-\lr s one may talk about embeddings. Isomorphisms and
automorphisms of \lr s are defined as expected (as those which respect the
ultrafilter $\cL$).
\end{Rmk}
\pa

\subsection{Confined \NNA, Separated and Rectified \Lr s}
\begin{Def}
Let $\cd{m}$ be an infinite cardinality. Given a \NA\ \AP, let $A$ be a set and $\U x\in\U A$. We say that $\U x$ is {\Dem $\cd{m}$-confined} or has {\Dem confinement $\cd{m}$} if there exists a (standard) set $A'\subset A$ of cardinality $\le\cd{m}$ such that $\U x\in\U A'$. If $\cd{m}$ is a cardinality such that all \NS\ $\U x$'s in $\U A$ for all sets $A$ are $\cd{m}$-confined, we say that the \NA\ is {\Dem $\cd{m}$-confined}, or that it has {\Dem confinement} $\cd{m}$.
\end{Def}

It is clear that for any infinite cardinality $\cd{m}$, and any (standard) map $f:A\to B$, $f$ maps \NS\ elements of confinement $\cd{m}$ to \NS\ elements of confinement $\cd{m}$ (if $\U x\in\U A'$ where $A'\subset A$ and $|A'|\le\cd{m}$, then $\U f(\U x)\in\U f(\U A')$. By Transfer $\U f(\U A')=\U (f(A'))$ and $|f(A')|\le\cd{m}$).\pa

Therefore if for each set $A$ we define $\U A_{\cd{m}}:=$ the set of $\U x\in\U A$ with confinement $\cd{m}$, we get a functor $A\mapsto\U A_{\cd{m}}$. That functor is a \NA. Indeed, for finite sets $F$ every \NS\ element is standard, thus $\U$belongs to a standard singleton, hence is of all confinements. Therefore $\U F_{\cd{m}}=\U F$; for a product $A\times B$, an $(\U x,\U y)\in\U (A\times B)=\U A\times\U B$ is $\cd{m}$-confined if and only if both $\U x$ and $\U y$ are so. Therefore
$\U (A\times B)_{\cd{m}}=\U A_{\cd{m}}\times\U B_{\cd{m}}$;
and for equalizers the required assertion is clear.\pa

We call the \NA\ $A\mapsto\U A_{\cd{m}}$: the \NA\ $A\mapsto\U A$
{\Dem $\cd{m}$-confined}, or the {\Dem $\cd{m}$-confinement} of the \NA\ $A\mapsto\U A$.\pa

Now suppose a \NA\ is given which has confinement $\cd{m}$, and let $B$ be some set of cardinality $\cd{m}$. Suppose we know the \lr\ $E=\U B$ over $B$ that $A$ induces on $B$.\pa

It seems evident that we can recover the \NS\ elements $\U A$ of any set (or class) $A$ using that \lr. Indeed, since the \NA\ has confinement $\cd{m}$, any \NS\ element $\U x$ of $A$ belongs to a set of the form $f(B)$ for some (standard) $f:B\to A$, hence, by Transfer, is of the form $\U x=f(e)$ for some $e\in E=\U B$. We can take the pair $(f,e)$ as a kind of ``coordinates'' for $\U x$, and all we need is to determine when two pairs $(f_1,e_1)$\,\, $(f_2,e_2)$ define the same $\U x$ and how $\U g:\U A\to\U A'$, for some $g:A\to A'$, is expressed using these ``coordinates''. But these characterisations are straightforward, since Transfer must hold. In fact, we use the ``template'' \lr\ $E$ over $B$ to define the \NS\ elements of any set by, in essence, ``carrying'' them over from $B$ in some analogy with the way tangent vectors may be defined in any smooth manifold by ``carrying'' them over from $\bR^n$.\pa

By this way we can start from any \lr\ $(E,\cL)$ over $B$ and define a $|B|$-confined \NA\, which will be (equivalent to) the {\em $|B|$-confinement} of any \NA\ which induces on $B$ a \lr\ isomorphic to $(E,\cL)$.\pa

Still, we can achieve that goal more ``neatly''.

\begin{Prop}\label{Prop:lrna}
Let $(E,\cL)$ be a \lr\ over an infinite set (basis) $B$. For any set $X$, let
$\U X$ be defined as the ``cylindrical'' ultrapower of $X$ with respect to $B^E$ and $\cL$, namely, the set $\left( X^{B^E}\right)_c$ of all mappings from $B^E$ to $X$ which depend only on a finite set of coordinates in $B^E$ (that set of coordinates depending on the mapping), factored by the equivalence relation of equivalence modulo $\cL$. (And the $\U$ of mappings between sets $f:X\to Y$ defined by composing with a representant of $\U x$).

Then: (i) $X\mapsto\U X$ is a $|B|$-confined \NA. It is called: the ($|B|$-confined) \NA\
{\Dem induced by the \lr}.

(ii) If $X\mapsto{}\Uz X$ is a $|B|$-confined \NA\ which induces on $B$ a \lr\ isomorphic to $(E,\cL)$, then $X\mapsto{}\Uz X$ is equivalent to $X\mapsto\U X$ via the two following natural mappings:

$\al=\al_X:\U X\to{}\Uz X$: let $\phi:B^E\to X$ be a representant of $\U x\in\U X$. $\phi$ depends only on a finite set of coordinates $e_1,\ldots,e_n$. Hence we may write for $\psi\in B^E$, $\phi(\psi)=f(\psi_{e_1},\ldots,\psi_{e_n})$,\,\,$f:B^n\to X$. We may assume $E={}\Uz B$, and define $\al(\U x):=\Big({}\Uz f\Big)(e_1,\ldots,e_n)\in{}\Uz X$.

$\ga=\ga_X:{}\Uz X\to\U X$: let $\Uz x\in{}\Uz X$. Since the \NA\ is $|B|$-confined, there is an $f:B\to X$ such that $\Uz x\in{}\Uz \Big(f(B)\Big)=\Big({}\Uz f\Big)(E)$. Thus $\Uz x=\Big({}\Uz f\Big)(e)$ for some $e\in E$. Define $\ga({}\Uz x)$ as $\psi\in B^E\mapsto f(\psi_e)$ modulo $\cL$. This does not depend on the choice of $f$.
\end{Prop}

\begin{Prf}
Clearly $X\mapsto\U X$ is a functor, and as in Example \ref{exm:na} one sees that it satisfies (i) and (ii) of Definition \ref{Def:na}. Therefore it is a \NA. Every $\U x\in\U X$, with a representant
$\phi:B^E\to X$ depending on the $n$ coordinates $e_1,\ldots,e_n$,
$\phi(\psi)=f(\psi_{e_1},\ldots,\psi_{e_n})$,
$\U$belongs to $f(B^n)$, of cardinality $\le|B|$. This proves (i).\pa

$\ga(\Uz x)$ does not depend on the choice of $f$. Indeed, if
$\Uz x=\Big({}\Uz f_1\Big)(e_1)=\Big({}\Uz f_2\Big)(e_2)$, then $(e_1,e_2)$\,\,$\Uz $belongs to the set $\{(b_1,b_2)\in B\times B\,|\,f_1(b_1)=f_2(b_2)\}$. This is the $\U$ of a $2$-relation in $B$, and by the way $\U$ is defined for a \lr, this means that for $\psi\in B^E$,\,\,$(\psi_{e_1},\psi_{e_2})$ is in that set modulo $\cL$, i.e.\ $f_1(\psi_{e_1})=f_2(\psi_{e_2})$ modulo $\cL$. This means that the two definitions for $\ga({}\Uz x)$ coincide.

Clearly $\al$ and $\ga$ are natural mappings.

$\al\circ\ga$ is the identity. Indeed, start from some $\Uz x\in{}\Uz X$.
$\U x:=\ga({}\Uz x)$ is $\psi\in B^E\mapsto f(\psi_e)$ modulo $\cL$, where
$f:B\to X$ is such that $\Uz x=\Big({}\Uz f\Big)(e)$. Thus $\psi\in B^E\mapsto f(\psi_e)$ is a representant of $\U x$, and $\al(\ga({}\Uz x))=\al(\U x):=\Big({}\Uz f\Big)(e)={}\Uz x$.

$\ga\circ\al$ is the identity. Indeed, start from some $\U x\in\U X$ represented by
$\phi:B^E\to X$ which depends on the coordinates $e_1,\ldots,e_n$, i.e.\
for $\psi\in B^E$,\,\,$\phi(\psi)=f(\psi_{e_1},\ldots,\psi_{e_n})$. Then
$\Uz x:=\al(\U x):=\Big({}\Uz f\Big)(e_1,\ldots,e_n)$. But by the $|B|$-confinement there are $g:B\to X$ and $e\in E$ such that
$\Uz x=\Big({}\Uz g\Big)(e)$. Then $(e,e_1,\ldots,e_n)$\,\,$\Uz $belongs to the set
$\{(b,b_1,\ldots,b_n)\in B^{n+1}\,|\,g(e)=f(e_1,\ldots,e_n)\}$, meaning that for
$\psi\in B^E$,\,\,$(\psi_e,\psi_{e_1},\ldots,\psi_{e_n})$ is in that set modulo $\cL$,
i.e.\ $g(\psi_e)=f(\psi_{e_1},\ldots,\psi_{e_n})$ modulo $\cL$. Thus, $\U x$ is represented also by $\psi\mapsto g(\psi_e)$, hence $\U x=\ga({}\Uz x)=\ga(\al(\U x))$.
\end{Prf}

\begin{Cor}\label{Cor:lrna}
Let $(E,\cL)$ be a \lr\ over an infinite $B$. Then TFAE: (i) There exists a \NA\ which induces on $B$ a \lr\ isomorphic to $E$. (ii) The \BC\ the \lr\ $E$ induces on $B$ a \lr\ isomorphic to $E$.

If these hold, we say that the \lr\ is {\Dem rectified}.

If these hold, one may take as an isomorphism the map $\ga_B:E\to\U B$ (notations as in Prop.\ \ref{Prop:lrna}).

For any $(E,\cL)$,\,\,$\ga_B$ can be characterized as mapping $e\in E$ to the projection $\psi\in B^E\mapsto\psi_e$ modulo $\cL$.
\end{Cor}

\begin{Prf}
If there exists a \NA\ $A\mapsto{}\Uz X$ inducing on $B$ a \lr\ isomorphic to $E$, then, by taking its $|B|$-confinement, we may assume it is $|B|$-confined, and then by Prop.\ \ref{Prop:lrna} it is equivalent to the $|B|$-confined \NA\ induced by $E$, hence the latter too induces on $B$ a \lr\ isomorphic to $E$.\pa

In this case we may take as the isomorphism the mapping $\ga_B$ from Prop.\ \ref{Prop:lrna}.\pa

For any $(E.\cL)$, in the definition of $\ga_B$ we identify ${}\Uz B=E$ and take
$f={\mathbf 1}_B$ (and $\Uz x=e$). Then $\ga_B(e)=(\psi\mapsto\psi_e)$ modulo $\cL$.
\end{Prf}

\begin{Def}
A \lr\ $E$ over $B$ is called {\Dem separated} if for $e_1,e_2\in E$, the set
$\{\psi\in B^E\,|\,\psi_{e_1}=\psi_{e_2}\}$ is in $\cL$ only if $e_1=e_2$.
\end{Def}

\begin{rmk}
For any \lr, the relation ``$\{\psi\in B^E\,|\,\psi_{e_1}=\psi_{e_2}\}$ is in $\cL$''
is an equivalence relation in $E$ and factoring $E$ by it gives a separated \lr.\pa
\end{rmk}

The \lr\ $E$ being separated just says that the map $\ga_B$ mentioned in Cor.\ \ref{Cor:lrna}: $e\mapsto\Big((\psi\mapsto\psi_e)\mod\cL\Big)$ is one-to-one. Therefore a rectified \lr\ is separated. (Or alternatively, the fact that the $\U$ of a diagonal is a diagonal means that the \lr\ that a \NA\ induces on a set $B$ must be separated.)\pa

For any \lr\ $E$ over $B$, the \lr\ of $B$ in the \BC\ it is called the {\Dem rectification} of $E$. Cor.\ \ref{Cor:lrna} can be phrased as saying that $E$ is rectified if and only if it is isomorphic to its rectification, and then as an isomorphism we can take the map
$\ga_B:E\to\U B$ mentioned there.\pa

\begin{Prop}\label{Prop:ga}
Let $(E,\cL)$ be a \lr\ over an infinite $B$, and suppose $E$ separated (then $\ga_B$ is injective. notations as above). Then $\ga_B$ is an isomorphism of \lr s over $B$ from $E$ to $\ga_B(E)\subset\U B$.
\end{Prop}

\begin{Prf}
Since $\ga_B$ is a bijection between $E$ and $\ga_B(E)$, all we need is to prove it preserves the ultrafilters on the cylinder Boolean algebras of $B^E$ and $B^{\ga_B(E)}$ which define the \lr s.\pa

Now, the set (here $h:B^n\to B$ is a function and $0\in B$ is some fixed (standard) element)
$$\Big\{\phi\in B^{\ga_B(E)}\,\Big|\,h\Big(\phi_{\ga_B(e_1)},\ldots,\phi_{\ga_B(e_n)}\Big)=0\Big\}$$
belongs to the ultrafilter of $\ga_B(E)$ if in $X\mapsto\U X=$ the \BC\ the \lr\ $E$, $$\U h\Big(\ga_B(e_1),\ldots,\ga_B(e_n)\Big)=0(=\nu(0)),$$
that is, using the characterisation of $\ga_B$ in Cor.\ \ref{Cor:lrna} (here $\psi$ varies over $B^E$),
$$\U h\Big(\Big(\psi\mapsto\psi_{e_1})\mod\cL\Big),\ldots,
\Big((\psi\mapsto\psi_{e_n})\mod\cL\Big)\Big)=(\psi\mapsto 0)\quad\mod\cL,$$
which means
$$\psi\mapsto h(\psi_{e_1},\ldots\psi_{e_1})=(\psi\mapsto 0)\quad\mod\cL,$$
i.e.,
$$h(\psi_{e_1},\ldots\psi_{e_1})=0\quad\mod\cL,$$
which just says that
$$\Big\{\psi\in B^E\,\Big|\,h(\psi_{e_1},\ldots,\psi_{e_n})=0\Big\}$$
belongs to $\cL$.
\end{Prf}

\subsection{Exact \Lr s}
\label{S:exact}
Abraham Robinson, the founder of \NA\ \cite{Robinson} did not, of course, wish that the \NA\ will be just the trivial $\U A=A$. His requirement that the \NA\ be an {\em enlargement} can be phrased, on the level of \lr s, as the requirement that for any ultrafilter $\cU$ on $B$ (or also on $B^I$,\,\,$I$ finite) there will be an element $\U x\in\U B$ (or in $(\U B)^I$) which is mapped to $\cU$ by the natural mapping that defines the \lr. But one may require slightly more:

\begin{Def}\label{Def:ex}
A \lr\ $E$ over an infinite $B$ is called {\Dem exact} if for any finite $I$, $J$ and
for any $\xi:I\to J$,\,\,the diagram expressing the naturality of the natural mapping between the (contravariant) functors $\bullet\to E^\bullet$ and $\be(B^\bullet)$ which defines the \lr, is ``exact'', in the sense that a member of $E^I$ and a member of $\be(B^J)$ which map to the same member of $\be(B^I)$ both come from some same member of $E^J$.\pa

Note, that for this to hold it suffices that it holds for inclusions ``adding one element'' $I\to I\cup\{i\}$ \NOT{guaranteeing that the transfer of relations commutes with
projections $B^{I\cup\{i\}}\to B^I$} and for the map $\{1,2\}\to\{1\}$
(the latter guarantees, as easily seen, that the $\U$ of the relation of equality in $B$ will be the relation of equality in $E$ -- i.e.\ that $E$ be separated).
\end{Def}

Note that for $I=\emptyset$ both $E^\emptyset$ and $\be(B^\emptyset)$ are singletons, $\{\emptyset\}$ if you wish, with subsets naturally viewed as the truth-values $\True=\{\emptyset\},\False=\emptyset$ which the natural mapping always preserves.\pa

If $E$ is exact then, in particular, applying the exactness of the diagram for the map $\emptyset\to\{1\}$ (or $\emptyset\to I$,\,\,$I$ finite) we have that for every ultrafilter $\cU$ in $B$ there is an element $e\in E$ whose ultrafilter, i.e.\ whose image in the natural mapping $E\to\be B$, is $\cU$, that is, $e$ $\U$belongs to all members of $\cU$ (and of course does not $\U$belong to non-members of $\cU$, which are the complements of members).%
\footnote{This does not hold, in general, for \NA\ constructed by an ultrapower
(if the power in the ultrapower is countable, a countable set will have only
$2^{\aleph_0}$ \NS\ members but has $2^{2^{\aleph_0}}$ ultrafilters.)}
If $\cU$ is fixed (principal) -- is the family of all subsets containing
some $a\in B$ -- than $e$ is unique and is the element of $E$ to be identified with $a$ (i.e.\ $e=\nu(a)$). If $\cU$ is free (non-principal), however, one proves that (if $E$ is exact) $e$ is never unique and we have ``conjugate'' \NS\ elements with the same ultrafilter which, as we shall see below, will be exchangeable by an automorphism with
a favorable choice of the \lr.
\pa

\begin{rmk}\label{rmk:Ii}
We may express the exactness of the diagram for a map $I\to I\cup\{i\}$ as follows:
For every family $\eta:I\to E$ of elements of $E$, and any ultrafilter $\cU$
in $B^{I\cup\{i\}}$ which projects on $B^I$ to the ultrafilter pulled-back from $\cL$ by $\eta$, one can find an element of $E$ as the image of $i$ so that the pull-back of $\cL$ by the extended $\eta$ will be $\cU$.\pa
\end{rmk}

This implies that for exact \lr s, $\U$ of relations will commute with projections $B^{I\cup\{i\}}\to B^I$ ($i\notin I$), i.e.\ with applying $\exists$. Of course, for any \lr, $\U$ of relations commutes with Propositional Calculus operations. Thus in an exact \lr\ it commutes with all first-order logical operations. And one has a {\em Transfer Principle}: Any true sentence made as a first-order logical combination of relations on $B$ will turn, by transferring each of the argument relations, into a sentence true for $E$.
\pa

\begin{Prop}
Any exact \lr\ $(E,\cL)$ over $B$ is rectified.
\end{Prop}
\begin{Prf}
By Prop.\ \ref{Prop:ga} all we need to prove is that if $E$ is exact the map $\ga_B$ mentioned there is onto.\pa

Recall that $\ga_B$ is a map from $E$ to $\U B$ = the set of cylindrical maps $\Psi:B^E\to B$ modulo $\cL$, which maps $e\in E$ to $\Big(\big(\psi\mapsto\psi_e\Big)\mod\cL\Big)\in\U B$.\pa

Now let $\U b\in\U B$ be some cylindrical map $\Psi:B^E\to B$ modulo $\cL$. Then
$$\Psi:\psi\in B^E\mapsto h\Big(\psi_{e_1},\ldots,\psi_{e_n}\Big)\in B,$$
where $e_1,\ldots,e_n\in E$ and $h$ is some $h:B^n\to B$.

Let $\cU$ be the push of $\cL$ by the map
$$\Phi_1:\psi\mapsto
\Big(\psi_{e_1},\ldots,\psi_{e_n},h\Big(\psi_{e_1},\ldots,\psi_{e_n}\Big)\Big)$$
to an ultrafilter in  $B^{n+1}$. Then $E$ being exact implies, by Remark \ref{rmk:Ii}, that there is an $e_0\in E$ such that the push of $\cL$ by the map
$$\Phi_0:\psi\mapsto\Big(\psi_{e_1},\ldots,\psi_{e_n},\psi_{e_0}\Big)$$
is the same ultrafilter $\cU$. Let
$$T:=\{(b_1,\ldots,b_n,b_0)\in B^n\,|\,b_0=h(b_1,\ldots,b_n)\}.$$
Then ${\Phi_1}^{-1}(T)$ is the whole $B^E$, hence $T\in\cU$, but that means that
${\Phi_0}^{-1}(T)\in\cL$, which just says that
$$\Psi(\psi)=h\Big(\psi_{e_1},\ldots,\psi_{e_n}\Big)=\psi_{e_0}\mod\cL.$$
Thus $\U b=\ga_B(e_0)$.
\end{Prf}

\subsection{$\cd{m}$-Exactness, $\cd{m}$-Saturation}
\begin{Def}\label{Def:mEx}
Let $\cd{m}$ is an infinite cardinal.

A \lr\ $E$ over an (infinite) set $B$ is called {\Dem $\cd{m}$-exact} if it is separated and exactness of the diagram of the natural mapping, as in Definition \ref{Def:ex}, holds for the inclusion $I\to I\cup\{i\}$ for any $|I|<\cd{m}$.

As in Remark \ref{rmk:Ii}, this means: For any $I$ such that $|I|<\cd{m}$, every family $\eta:I\to E$ of elements of $E$, and any ultrafilter $\cU$ in the cylinders of $B^{I\cup\{i\}}$ which projects on $B^I$ to the ultrafilter pulled-back from $\cL$ by $\eta$, one can find an element of $E$ as the image of $i$ so that the pull-back of $\cL$ by the extended $\eta$ will be $\cU$.
\end{Def}

\begin{rmk}\label{rmk:M}
By a transfinite process of adding elements $i$, the $\cd{m}$-exactness of $E$ will
guarantee exactness of the diagram of the natural mappings for any inclusion $I\to M$ for $|I|<\cd{m},|M|\le\cd{m}$.
\end{rmk}

In particular, by the definition, a \lr\ is exact if and only if it is separated and $\aleph_0$-exact.\pa

\begin{Prop}\label{Prop:BmExB}
Let $\cd{m}$ is an infinite cardinal.

Let $E$ be an $\cd{m}$-exact \lr\ over an (infinite) set $B$. Let $B'\subset B$,\,$B'$ also infinite. Let $E'=\U B'$ be the set of members of $E$ that $\U$belong to $B'$. Then the sub-\lr\ $E'$ over $B'$ is also $\cd{m}$-exact.
\end{Prop}

\begin{Prf}
Let, as in Def.\ \ref{Def:mEx}, $I$ satisfy $|I|<\cd{m}$,\,\,$\eta:I\to E'$ and $\cV$ an ultrafilter in the cylinders of $(B')^{I\cup\{i\}}$ which projects on $(B')^I$ to the ultrafilter pulled-back by $\eta$ from the defining cylinder ultrafilter $\cL'$ of $E'$.

Note, that if $\cL$ was the defining cylinder ultrafilter of $E$, then for any finite set $e_1,\ldots,e_n\in E'$, the preimage of $(B')^n\subset B^n$ by
$B^I\to B^{\{e_1,\ldots,e_n\}}$ is a member of $\cL$ and $\cL'$ obtains from $\cL$ by intersecting its members with these preimages. Thus, the inclusion $(B')^{I\cup\{i\}}\to B^{I\cup\{i\}}$ maps $\cV$ into an ultrafilter in the cylinders of $B^{I\cup\{i\}}$ (which contains all primages of $(B')^n$ by $B^I\to B^{\{\iota_1,\ldots,\iota_n\}}$ for
$\iota_1,\ldots,\iota_n\in I\cup\{i\}$) and which will project on $B^I$ to the
ultrafilter pulled-back by $\eta$ from $\cL$.

Hence, by the $\cd{m}$-exactness of $E$, there is an element of $E$ as $\eta(i)$ so that the pull-back of $\cL$ by the extended $\eta$ is $\cU$. Since $\cU$ contains the preimage of $B'$ by the projection $B^{I\cup\{i\}}\to B^{\{i\}}$, $\eta(0)$ belongs to $\U B'=E'$.
Thus the extended $\eta$ maps $I\cup\{i\}$ to $E'$. And the way $\cU$ is constructed from $\cV$ and $\cL'$ is constructed from $\cL$ means that the pull-back of $\cL'$ by the extended $\eta$ is $\cV$, as required.
\end{Prf}

$\cd{m}$-exactness is related to the property of a \NA\ being {\em $\cd{m}$-saturated}.\pa

Recall (see end of \S\ref{S:NAFco}) the distinction between {\Dem internal} and {\Dem external} sets of \NS\ elements: For some (standard) set $A$, a subset of $\U A$ is called internal if it is of the form $\HAT{\U X}:=\{\U y\in\U A\,|\,\U y\,\,\T\in\U X\}$ for some $\U X\in\U\cP A$ where $\cP A$ is the power set. A subset of $\U A$ which is not internal is called external.\pa

\begin{rmk}
An alternative characterization for an internal set: an internal subset of $\U A$ is any set of the form: the section parallel to $\U A$,\,$\U R[\U y]$ of $\U R$ through some element $\U y\in\U C$ for some (standard) set $C$ and some (standard) relation $R\subset C\times A$.\pa

Indeed, by Transfer, $\U R[\U y]$ is just the set of \NS\ $\U$members of the image $\U X$ of $\U y$ by the $\U$ of the map $y\in C\mapsto R[y]\in\cP A$. On the other hand, if $\U X\in\U\cP A$, then the set of its $\U$members is the section parallel to $\U A$ through $\U X$ of the $\U$ of the relation $y\in X$ in $A\times\cP A$.
\end{rmk}

\begin{Def}
Let $\cd{m}$ is an infinite cardinal.

A \NA\ is called {\Dem $\cd{m}$-saturated inside a (standard) set $A$}, if a family of cardinality $<\cd{m}$ of {\em internal} subsets of $\U A$ has non-empty intersection provided any finite subfamily has non-empty intersection. If a \NA\ is $\cd{m}$-saturated inside every (standard) set, then it is called {\Dem $\cd{m}$-saturated}.
\end{Def}

\begin{Prop}\label{Prop:ExSat}
Let $\cd{m}$ is an infinite cardinal.\pa

(i) Suppose a \lr\ $E$ over an infinite $B$ is $\cd{m}$-exact. Then the \BC\ it is $\min(\cd{m},|B|^+)$-saturated, while it is $\cd{m}$-saturated inside every (standard) set of cardinality $\le|B|$.\pa

(ii) If a \NA, \AP, is $\cd{m}$-saturated, and $B$ is an infinite (standard) set, then the \lr\ which the \NA\ induces over $B$ is $\cd{m}_1$-exact if $\cd{m}'<\cd{m}_1\Rightarrow 2^{|B|}\cdot\cd{m}'<\cd{m}$.
\end{Prop}

\begin{Prf}
(i) Note first, that since the \lr\ $E$ is $\cd{m}$-exact, hence exact, it is rectified. Therefore we may identify $\U B=E$.

Let us be given a family of cardinality $<\cd{m}$ of internal subsets of the $\U$ of some standard set $C$, where every finite subfamily has non-empty intersection.

Every member of the family of internal sets is the scope $\HAT{\U X_\iota}$ of some $\U X_\iota\in\U\cP C$, and, the \NA\ being $|B|$-confined, there is a standard subset $Y_\iota\subset\cP C$ of cardinality $\le|B|$ such that $\U X\in\U Y_\iota$. Denoting by $R_\iota$ the relation $\in$ in $C\times Y_\iota$, $\HAT{\U X_\iota}$ is the section parallel to $\U C$ of $\U R_\iota$ through $\U X_\iota$. As $|Y_\iota|\le|B|$, we may
replace the $Y_\iota$'s with subsets of $B$, shifting the $R_\iota$'s accordingly.
Then the $\U X_\iota$'s belong to $E=\U B$, and they form an element of $E^I$ for the
set $I$ of the $\iota$'s, where $|I|<\cd{m}$.

Now we distinguish between the two cases. If $|C|\le|B|$, we may assume $C\subset B$.
If we do not have $|C|\le|B|$, replace $\cd{m}$ by $\min(\cd{m},|B|^+)$, that is, assume
$\cd{m}\le|B|^+$. Then the cardinality of the family of internal sets is $\le|B|$ and if we pick an element from every intersection of a finite subfamily, we get at most $|B|$ elements which, by $|B|$-confinement, all belong to some $\U C'$ for some (standard) $C'\subset C$ with $|C'|\le|B|$. And we proceed with $\HAT{\U X_\iota}\cap\U C'$ instead of $\HAT{\U X_\iota}$ as in the first case.

So we assume $C\subset B$; writing $B=B^{\{0\}}$, every $R_\iota$ is a relation $=$ subset of $B^{\{0,\iota\}}$. Which can be viewed as a ``cylinder'' in $B^{\{0\}\cup I}$. We wish to have an ultrafilter $\cU$ in the cylinders of $B^{\{0\}\cup I}$ which contains these $R_\iota$'s and also all sets which depend only on a finite number of coordinates {\em in $I$}, that the $\U X_\iota$ satisfy. The fact that $\cU$ exists, i.e.\ that any finite number of these intersect, follows from the fact that their $\U$ intersect because any finite subfamily of the $\U R_\iota[\U X_\iota]$ intersects. Now, $\cd{m}$-exactness finishes the task: it says that there is a $\U X_0\in E$ such that the push of $\cL$ to $B^{\{0\}\cup I}$ is $\cU$, that is, the family $(\U X_0,(\U X_\iota))$ $\U$belongs to all members of $\cU$, meaning that $\U X_0$ is in the intersection of all the members of the family of internal sets.\pa

(ii) Denote $E=\U B$. Let $I$ be such that $|I|=\cd{m}'<\cd{m}_1$, let $\eta:I\to E$ be a family of elements of $E$, and $\cU$ an ultrafilter in the cylinders of $B^{I\cup\{0\}}$ which projects on $B^I$ to the ultrafilter $\cU'$ pulled-back from $\cL$ by $\eta$, and one wishes to find an element of $E$ as $\eta(0)$ -- the image of $0$, so that the pull-back of $\cL$ by the extended $\eta$ will be $\cU$.

To satisfy this, $\eta(0)$, together with the $\eta(\iota),\iota\in I$, must $\U$belong to all the members $U$ of $\cU$. Belonging to each such $U$ requires $\eta(0)$ to belong to a certain internal subset of $\U B=E$. To belong to any finite number of these internal subsets, is to satisfy the conjunction of these $\U U$'s with a finite number of $\eta(\iota)$'s, $\eta_{\iota_1},\ldots,\eta_{\iota_n}$ as parameters. The set of these parameters belongs to the $\U$ of a any member $U'$ of $\cU'$. Hence, Transfer will tell us that there is such an $\eta(0)$, if we can find a member $U'$ of $\cU'$ so that for any ``standard'' ``parameters'' $b_{\iota_1},\ldots,b_{\iota_n}\,\in B$ there is a $b_{\iota_0}\in B$ which together with them belong to all the (finite number of) $U$'s. But for this just take as $U'$ the projection on $B^I$ of the intersection of these $U$'s.

Thus, every finite number of these internal subsets intersect.

Also, the cardinality of $\cU$ is at most the cardinality of the set of all ``cylinders'' in $B^I$, which is $\le 2^{|B|}\cdot\cd{m}'<\cd{m}$.

Hence, by saturation, all requirements for $\eta(0)$ can be fulfilled.
\end{Prf}

\subsection{Regular \Lr s, \Lr s with Homogeneity and Universality}
\begin{Def}
A \lr\ $E$ over $B$ is called {\Dem regular} if it is $|E|$-exact.
\end{Def}

A regular \lr\ is {\em universal}, in the sense that any embedding in $E$ of a sub-\lr\ $M'$ of cardinality $<|E|$ of a separated \lr\ $M$ over $B$ of cardinality $|M|\le|E|$ can be extended to an embedding of all $M$.\pa

Indeed, the embedding of $M'$ can be viewed as a family of elements of $E$ so that the pull-back of the defining cylinder ulterafilter in $B^E$ is the defining cylinder ultrafilter in $B^{M'}$, which is the projection of the defining cylinder ultrafilter in $B^M$, and one uses Remark \ref{rmk:M}. (Separation will guarantee that distinct elements map to distinct elements, so one really has a one-to-one embedding)\pa

If, instead of the argument of Remark \ref{rmk:M} one uses a transfinite back-and-forth construction, i.e.\ adds an element $i$ on one side of the embedding, then an element on the other, then again on the first side etc., one concludes that a regular \lr\ is {\em homogeneous} -- any isomorphism between two sub-\lr s of cardinality $<|E|$ in two regular \lr s $E$ over $B$ with the same cardinality $|E|$, can be extended to an isomorphism between the $E$'s. In particular (if we start from the empty isomorphism) any two such \lr s $E$ are isomorphic. Also (starting from singletons) any two members of $E$ which map to the same ultrafilter in $\be B$ are exchangeable by an automorphism of $E$. \pa

Thus one might say that there is, up to isomorphism, at most one regular \lr\ over $B$ of any specific cardinality.\pa

Prop.\ \ref{Prop:BmExB} implies

\begin{Prop}\label{Prop:BregB}
Let $E$ be a regular \lr\ over an (infinite) set $B$. Let $B'\subset B$,\,$B'$ also infinite. Let $E'=\U B'$ be the set of members of $E$ that $\U$belong to $B'$. Then the sub-\lr\ $E'$ over $B'$ is also regular.
\end{Prop}

Indeed, by Prop.\ \ref{Prop:BmExB} $E'$ is $|E|$-exact, {\em a fortiori} $|E'|$-exact.

\subsection{Cardinalities}
\label{S:car}
\begin{Prop}
Suppose a \NA\ has confinement $\cd{m}$ and saturation $\cd{m}'$. Then $\cd{m}'\le\cd{m}^+$.
\end{Prop}

\begin{Prf}
Take a (standard) set $A$ of cardinality $\cd{m}^+$, and well-order it by the set of all ordinals of cardinality $\le\cd{m}$. Take the standard subsets $[a,\to[$ of $A$, and the corresponding internal sets of their $\U$members. Any finite number of these intersect, and there are $\cd{m}^+$ of them. If we had saturation $(\cd{m}^+)^+$ there would be a $\U$member $\U a$ of $A$ which belongs to all of them, hence is bigger than all standard
$a\in A$, and thus cannot $\U$belong to any standard set $A'\subset A$ of cardinality $\le\cd{m}$ since any such set is bounded by some (standard) $a\in A$. Thus we would not have confinement $\cd{m}$.
\end{Prf}

If a \lr\ $E$ over $B$ is exact, then, as we saw in \S\ref{S:exact}, for every ultrafilter $\cU$ in $B$ there is a $\U b\in E$ that maps to $\cU$ in $E\to\be B$. Since there are
$2^{2^{|B|}}$ ultrafilters in $B$ (see Appendix \ref{A:U}),  we have

\begin{Prop}\label{Prop:ExactCar}
If a \lr\ $E$ over $B$ is exact, then $|E|\ge 2^{2^{|B|}}$.
\qed\end{Prop}

\begin{Prop}\label{Prop:22b}
Let $B$ be an infinite set of cardinality $\cd{b}$. Then there exists a \lr\ $E$ over $B$ such that:
(i) $|E|=2^{2^\cd{b}}$; (ii) The \lr\ $E$ is $(2^\cd{b})^+$-exact.

(This does not require any Generalized Continuum Hypothesis (GCH). But if we have GCH, at least the case $2^{2^\cd{b}}=(2^\cd{b})^+$ of GCH, we conclude that $E$ is regular -- is {\em the} regular \lr\ over $B$ with cardinality $2^{2^\cd{b}}$).

By Prop.\ \ref{Prop:ExSat} we will have that the \BC\ this \lr\ is $\cd{b}^+$-saturated, and is $(2^\cd{b})^+$-saturated inside a (standard) set of cardinality $\le\cd{b}$.
\end{Prop}

\begin{Prf}
Let $\Omega$ be the set of ordinals of cardinality $\le2^\cd{b}$. It is the smallest ordinal of cardinality $(2^\cd{b})^+$.

We shall construct $E$, together with the defining cylinder ultrafilter $\cL$, in steps indexed by $\Omega$, always adding $e\in E$'s so that for a family $\eta:I\subset 2^B\to$ the set $E'$ of previously constructed $e$'s and for a cylinder ultrafilter $\cU$ in $B^{I\cup\{0\}}$ whose projection on $B^I$ is the ultrafilter pulled-back from $\cL$ by $\eta$, the added $e$ will serve as $\eta(0)$ so that $\cU$ is the ultrafilter pulled-back by the extended $\eta$. To this end we appropriately choose the ultrafilter $\cL$ in the cylinders of $B^{E'\cup\{e\}}$, of course one that projects on $B^{E'}$ to the previously constructed $\cL$.

For any $\iota\in\Omega$, supposing we already have $E_\iota$ and $\cL_\iota$, we shall take care simultaneously of all possible $\eta$ into $E'=E_\iota$ and $\cU$'s, adding one $e$ for each of them, to get $E_{\iota+1}$. If we assume as an induction hypothesis that $|E_\iota|\le2^{2^\cd{b}}$, then there are at most
$\left( 2^{2^\cd{b}}\right)^{2^\cd{b}}=2^{2^\cd{b}}$ such $\eta$'s, and for each $I\subset 2^B$ the cardinality of the cylinder Boolean algebra of
$B^I$ is $|I|\cdot2^\cd{b}\le2^\cd{b}$, hence the cardinality of all possible cylinder ultrafilters, these being subsets of this Boolean algebra, is $\le2^{2^\cd{b}}$. Thus there are at most $2^{2^\cd{b}}$ pairs $\eta,\cU$ to take care of, hence $2^{2^\cd{b}}$ added $e$'s, and also $|E_{\iota+1}|\le2^{2^\cd{b}}$.

For limit ordinals $\iota<\Omega$ (and, in fact, also for $\Omega$ itself), we just take as $E_\iota$ the union of everything constructed at $\iota'<\iota$ (a union of at most $2^\cd{b}$ terms -- $(2^\cd{b})^+$ for $\Omega$) and as $\cL$ the union of the $\cL$'s -- note that we are always considering cylinders, which depend only on a finite number of coordinates. Clearly we will preserve the fact $|E_\iota|\le2^{2^\cd{b}}$.

The final $E$ will be $E_\Omega$, after we had made it {\em separated} by taking its factor set with respect to the equivalence relation on $e_1,e_2$: $\{\psi\,|\,\psi_{e_1}=\psi_{e_2}\}\in\cL$.
(i) will be satisfied. (By the above $|E|\le2^{2^\cd{b}}$. the inverse inequality is Prop.\ \ref{Prop:ExactCar}).
(ii) will be satisfied too, since every $\eta$ for $E_\Omega$ concerns at most $2^\cd{b}$ elements $e$, constructed at $2^\cd{b}$ $\iota$'s, the set of these $\iota$'s is bounded in $\Omega$ by some $\iota'$, and this $\eta$ has been dealt with -- received its $\eta(0)$ -- already in $E_{\iota'}$.

We still have to prove that we can construct the cylinder ultrafilter $\cL_{\iota+1}$ when we construct $E_{\iota+1}$ above. $\cL_{\iota+1}$ must: (a) Project to $\cL_\iota$ on
$B^{E_\iota}$. That is, contain all cylinders, depending on a finite number of coordinates in $E_\iota$, that belong to $\cL_\iota$. (b) For each $e$ constructed in the $\iota+1$ step for some $\eta$ and $\cU$, it must be pulled back by the extended $\eta$ to $\cU$. That is, for every $U\in\cU$ which depends, as a subset of $B^{n+1}$, on the coordinates $0,i_1,\ldots,i_n$, with $\eta(i_j)=e_j\in E_\iota$, $j=1,\ldots,n$ it must contain the set
\begin{equation}\label{eq:psiOmega1}
\{\psi\in B^{E_{\iota+1}}\,|\,(\psi_e,\psi_{e_1},\ldots,\psi_{e_n})\in U\}.
\end{equation}
To prove such a cylinder ultrafilter $\cL_{\iota+1}$ exists, all we need to show is that any finite set of these sets, which it must contain, intersect. That is, we must prove that for any cylinder $V$ in $\cL_\iota$, which is of the form (here $e_1,\ldots,e_n\in E_\iota$, $V\subset B^n$):
\begin{equation}\label{eq:psiOmega2}
\{\psi\in B^{E_{\iota+1}}\,|\,(\psi_{e_1},\ldots,\psi_{e_n})\in V\},
\end{equation}
and any finite number of $e$'s and sets as in (\ref{eq:psiOmega1}), there is a $\psi$ which belongs to all of them. For the latter, we will be able to find a $\psi_e$ if $(\psi_{e_1},\ldots,\psi_{e_n})$ belongs to the projection $U'$ of $U$ on $B^n=B^{\{i_1,\ldots,i_n\}}$. This projection is, by assumption, in the cylinder ultrafilter pulled-back from $\cL_\iota$ by the unextended $\eta$. This means that the set of $\psi$ satisfying this is a cylinder in $\cL_\iota$, as in (\ref{eq:psiOmega2}). Hence we are left with a finite number of conditions of the form (\ref{eq:psiOmega2}), for cylinders in $\cL_\iota$ which intersect, $\cL_\iota$ being a cylinder ultrafilter.
\end{Prf}

\begin{Prop}\label{Prop:car}
Let $\cd{b}$ be an infinite cardinality.

Let us work in a \NA\ with confinement $\cd{b}$, saturation $\cd{b}^+$ and saturation
$(2^\cd{b})^+$ inside (standard) sets of cardinality $\le\cd{b}$, and in which for $|B|=\cd{b}$ we have $|{}\U B|=2^{2^\cd{b}}$. Then
\begin{itemize}
\item[(i)]
Let $B$ be a set of cardinality $\cd{b}$. Then for every internal subset $\mathbf\be\subset\U B$ which contains all standard elements of $B$, and for every ultrafilter $\cU$ on $B$, $\exists\;\U x\in\mathbf\be$ whose ultrafilter is $\cU$.

Consequently, $|\mathbf\be|=2^{2^\cd{b}}$ (recall that $|\U{}B|=2^{2^\cd{b}}$).
\item[(ii)]
If $\mathbf\al$ is internal, then either $\mathbf\al$ is finite or $\exists$ an internal one-one mapping from an internal $\mathbf\be$ as in (i) to $\al$, hence $|\mathbf\al|\ge2^{2^\cd{b}}$.
\item[(iii)]
If $\mathbf\al$ is $\U$finite (i.e.\ a $\U$member of a standard family of finite
sets) then $|\mathbf\al|$ is finite or $2^{2^\cd{b}}$.
\end{itemize}
\end{Prop}

\begin{Prf}
(i) We use $(2^\cd{b})^+$-saturation (which holds when working inside $B$):
for any finite family of members of $\cU$\,\,$\exists\;\U x\in\mathbf\be$ belonging
to all of them (indeed a standard one), and by saturation we are done.

(ii) Here we use $\cd{b}^+$-saturation: If $\mathbf\al$ is not finite, then for every
finite set $F$ of standard members of $B$\,\,\,$\exists$ an internal set
$\mathbf\be_F\subset\U B$ containing all members of $F$ and an internal one-one mapping
$\mathbf\be_F\to\al$ (take $\mathbf\be_F=F$), hence by saturation $\exists$ a $\mathbf\be$ suitable for all standard members of $B$.

(iii) Suppose $\mathbf\al$ not finite. By (ii) $|\mathbf\al|\ge2^{2^\cd{b}}$. On the other hand, by $\cd{b}^+$-confinement $\mathbf\al$ $\U$belongs to a standard set $S$ of cardinality $\le\cd{b}$ whose members are finite sets. Then $|\cup S|\le\cd{b}$, thus
$|{}\U(\cup S)|\le2^{2^\cd{b}}$, while $\mathbf\al\subset{}\U(\cup S)$.
\end{Prf}

From Prop.\ \ref{Prop:22b}, Prop.\ \ref{Prop:car} and Prop.\ \ref{Prop:BregB} we get

\begin{Cor}
Let $E$ be a \lr\ over $B$ as in Prop.\ \ref{Prop:22b}. Let $B'$ be an infinite subset of $B$. Suppose we have the case $2^{2^{|B|}}=\Big(2^{|B|}\Big)^+$ of GCH. Let $E'=\U B'\subset E$. Then the \lr\ $E'$ over $B'$ is {\em the} regular \lr\ of cardinality $2^{2^{|B|}}$ over $B'$.
\end{Cor}

\section{Iterated \NA}
\subsection{The Embeddings $\nu$ and $\U\nu$}
Fix a \NA\ \AP. For every set $A$, we have $\U A$, but, of course, $\U A$ being again a set, we have $\U\U A$, also $\U\U\U A$ etc. The \NA\ defines, for every set $A$, an injection $\nu:A\to\U A$ mapping $a\in A\mapsto\nu(a)\in\U A$. Therefore we have (again an injection, by Transfer) $\U\nu:\U A\to\U\U A$. Note that we still have $\nu:\U A\to\U\U A$.\pa

Both $\nu$ and $\U\nu$ are injective. Thus both ``embed'' $\U A$ in $\U\U A$. Yet they are highly different, as we shall see.\pa

$\nu$, of course, satisfies Transfer. But by Transferring Transfer we find that also $\U\nu$ satisfies a kind of Transfer. One just has to be careful with quantifiers -- on which set to take them.\pa

Let us note that
\begin{Prop}
For any (standard) set $A$, and any (standard) $a\in A$, the $\U\U$A-elements $\nu(\nu(a))$ and $(\U\nu)(\nu(a))$ are the same.
\end{Prop}
\begin{Prf}
For any map $f:A\to B$, $\U f:\U A\to\U B$ does to the $\nu$ of elements of $A$ the same as $f$ does to the elements. This means $\U f(\nu(a))=\nu(f(a)),\,a\in A$. Now take $B=\U A$ and $f:A\to B$ to be $\nu:A\to\U A$.
\end{Prf}
And one easily sees that $a\mapsto\nu(\nu(a))$ and $a\mapsto(\U\nu)(\nu(a))$ ``transfer'' relations and functions in the same way.\pa

To clarify the relationship between these $\nu$ and $\U\nu$ let us consider, as an example, the case $A=\bR$, with its subset $\bN$.\pa

\subsection{\NNA\ and Iterated \NNA\ in $\bR$ and $\bN$}
We place ourself in a \NA\ which is $\aleph_0^+$-saturated, and look at the \NA\ inside $\bR$, in particular inside its subset $\bN$.\pa

$\aleph_0^+$-saturation implies that for a descending sequence $\U a_n$, (and an ascending sequence $\U b_n$) of $\U$members of $\bR$ (or any other linearly ordered set), such that $\U a_n>\U b_m\,\forall n,m$  there is a $\U$member $\U c$ of $\bR$ such that $\U a_n>\U c$ (and $\U c>\U b_m$)\,$\forall n,m$. In particular, there are $\U$reals, also $\U$natural numbers, that are bigger than all finite standard numbers. These are sometimes called (positive) {\em infinite} or {\em unlimited}, while a $\U$real whose absolute value is not unlimited will be called {\em limited}. While for any $n\in\bN$ the set $]\leftarrow,n]$ in $\bN$ is finite, therefore all its $\U$membrs are standard, and consequently all $\U$natural numbers are either standard or unlimited, there are $\U$real numbers $>0$ but smaller than any standard positive real. As is well-known, these are called (positive) {\em infinitesimal} $\U$reals (and a positive or negative $\U$real is called infinitesimal if its absolute value is infinitesimal). As is well-known, and easily proved using the completeness of $\bR$ (say, the existence of suprema and infima), for any {\em limited} $\U a\in\U\bR$ there is a unique standard $\textrm{st}(\U a)\in\bR$ so that $\U a-\textrm{st}(\U a)$ is infinitesimal. Two $\U$reals whose difference is infinitesimal are called {\em near}.\pa

Now, as said above, $\U\bR$ also has its $\U$members, forming the set $\U\U\bR$ which may be referred to as the $\U\U$reals. And there is a $\nu:\U\bR\to\U\U\bR$ which maps each $\U$real to the $\U\U$real to be identified with it, and since there is $\nu:\bR\to\U\bR$, it gives rise to $\U\nu:\U\bR\to\U\U\bR$. As said above, these are two
``embeddings'' of $\U\bR$ into $\U\U\bR$.\pa

Let's take a look at the set of positive (thus nonzero) infinitesimal $\U$reals. Note that, by the above, for any ascending (resp.\ descending) sequence of them there is a positive infinitesimal further bigger (resp. smaller) than all of them. The set of infinitesimals is external -- if it were internal it would have a $\U$supremum which is impossible. It has $\U$members, the $\U$infinitesimals -- thay are $\U\U$reals. In particular, since for any set $A$\,\,$\nu$ maps $A$ into $\U A$, $\nu$ of infinitesimals are $\U$infinitesimals.\pa

Now, since the positive infinitesimals are the positive $\U$reals that are smaller than $\nu(r)$ for all reals $r>0$, we have by Transfer: the positive $\U$infinitesimals are the positive $\U\U$reals smaller than $\U\nu(\U r)$ for all $\U$reals $\U r>0$. But note that we can take here as $\U r$ a positive infinitesimal! Thus we have:\pa

All positive $\U$infinitesimals, which include the $\nu$ of positive infinitesimals, are smaller than all $\U\nu$ of positive infinitesimals!\pa

Turning to the other end, there are the positive $\U$unlimited $\U\U$reals, which include the $\nu$ of unlimited $\U$reals. They are characterized as bigger than $\U\nu(\U r)$ for all $\U r\in\U\bR$. Thus they are bigger than the $\U\nu$ of unlimited $\U$reals.\pa

Consider, in particular, $\U\U$natural numbers. We saw that the $\U$naural numbers are partitioned into two sets: $\nu(\bN)$ and the unlimited $\U$natural numbers, the latter all bigger than the former. This gives by Transfer: the $\U\U$natural numbers are partitioned into $\U\nu(\U\bN)$ and the $\U$unlimited $\U\U$natural numbers, the latter all bigger than the former. The latter include the $\nu$ of unlimited $\U$natural numbers. The former are, of course, partitioned into $\nu(\nu(\bN))=(\U\nu)(\nu(\bN))$, and the $\U\nu$ of the unlimited, the latter all bigger than the former. To repeat: $\U\U\bN$ is partitioned into three parts: $\nu(\nu(\bN))=(\U\nu)(\nu(\bN))$, all smaller than $\U\nu$ of the unlimited, these all smaller than the $\U$unlimited, which include the $\nu$ of the unlimited.\pa

If we assume the \NA\ is $\aleph_0$-confined (or more generally, is $\cd{m}$-confined and $\cd{m}^+$-saturated for some infinite cardinal $\cd{m}$), we will have: every positive $\U$infinitesimal $\U\U\eps$\,$\U$belongs to a countable set of infinitesimals, the latter, by $\aleph_0^+$-saturation, as said above, has a positive infinitesimal upper (resp.\ lower) bound. Hence (the $\nu$ value) of the bound is bigger (resp.\ smaller) also than $\U\U\eps$. Thus ``being less than (the $\nu$-value of) a positive infinitesimal'' is the same as ``being less than a positive $\U$infinitesimal''. Similarly, in $\U\U\bR$ or in $\U\U\bN$, every $\U$unlimited $\U\U$number has an upper (resp.\ lower) bound which is $\nu$ of an unlimited $\U$number.\pa

We had the (equivalence) relation in $\U\bR$ of ``being near'', i.e.\ having an infinitesimal difference. In fact, it is an intersection $\cap_n\U U_n$ of the $\U$'s of a decreasing sequence of standard relations $U_n\subset\bR\times\bR$, in this case a basis to the uniformity in $\bR$. We had: the equivalence class of a limited $\U a\in\U\bR$ contains a unique standard $\nu(\textrm{st}(\U a))$,\,\,$\textrm{st}(\U a)$ called the {\em standard part} of $\U a$. Moving from such a relation to a relation in $\U\U\bR$, we have two options:\pa

The ``strong'' (equivalence) relation $\U(\cap_n\U U_n))$ -- the difference being a $\U$infinitesimal. The existence of the standard part implies by Transfer: the equivalence class of a $\U$limited $\U\U a\in\U\U\bR$ contains a unique element $\U\nu(\U\textrm{st}(\U\U a))\in\U\nu(\U\bR)$.\pa

The ``weak'' (equivalence) relation $\cap_n\U\U U_n)$, in our case saying that the difference being less than any $\nu(\nu(a))$,\,\,$a>0,a\in\bR$. Lets call such $\U\U$reals {\em ``weakly infinitesimal''}. These are the elements equivalent to $0$, and they include all $\U\nu$ of infinitesimals. In fact, any $\U\U a\in\U\U\bR$ with standard bounds is thus ``weakly'' equivalent to $\nu(\nu\big(\textrm{st}(\U\textrm{st}(\U\U a))\big))$, which is the only standard element in its equivalence class.\pa

(Note that we could start with only the rational numbers $\bQ$, consider $\U\bQ$, define infinitesimals and being near as above, and then recover $\bR$ as the set of equivalence classes of {\em limited} $\U$rationals with respect to nearness. Transferring, $\U\bR$ will be identified with the set of equivalence classes of $\U$limited $\U\U$rationals with respect to $\U$nearness, i.e.\ having a $\U$infinitesimal difference.)

\subsection{A Criterion for $\nu$ and $\U\nu$ of $\U$Elements to Satisfy Some Standard Relation, and Applications}

Return to general (standard) sets.
\begin{Prop}\label{Prop:nunustar}
Let $A$ and $B$ be (standard) sets and let $R\subset A\times B$ be a (standard) relation between elements of $A$ and elements of $B$. Let $\U a\in\U A$ and $\U b\in\U B$. Then a necessary and sufficient condition that $\U\nu(\U b)$ is in the relation $R$ with $\nu(\U a)$ is
$$\U b\,\,\U\textrm{belongs to the set of (standard) }b\in B\textrm{ that are in the relation }R\textrm{ with }\U a\eqno{(\star)}$$
\end{Prop}
\begin{Prf}
Indeed, $(\star)$ holds if and only if $\U b$ satisfies the $\U$ of the property of $b\in B$: ``$\nu(b)$ is in the relation $R$ with $\U a$''. By Transferring, this means that it satisfies: ``$\U\nu(\U b)$ is in the relation $R$ with $\nu(\U a)$''.
\end{Prf}
{\EM Some examples}:
\begin{itemize}
\item Take $R$ as the relation of {\em inequality} in $A$. Thus to have $\U\nu(\U b)\ne\nu(\U a)$,\,\,$\U b$ must $\U$beloing to the set of $b\in A$ unequal to $\U a$, this set being the whole $A$ unless $\U a$ is standard. For standard $\U a=\nu(a),\,a\in A$,\,\,$\U\nu(\U b)$ can equal $\nu(\U a)=\U\nu(\nu(a))$ only if $\U b=\nu(a)$. Therefore:\pa

The $\U\nu$ of a $\U$A-member can equal the $\nu$ of a $\U$A-member only if both are standard and equal.\pa

\item Let $A$ be partially ordered, and take as $R$ the relation $a\ge b$. $(\star)$ says that $\U b$ $\U$belongs to the cone of the standard minorants of $a$. Thus if $\U a$ is ``infinite'' in the sense of majorizing all standard elements, then the latter cone is everything so we always have $\U\nu(\U b)\le\nu(\U a)$ (similarly to what we saw in $\bR$ and $\bN$ in the previous subsection).\pa

\item Let $A$ be a (standard) set and $B$ a set of subsets of $A$, with the relation $a\in b$. Then $\nu(\U a)\in\U\nu(\U b)$ (here we identify a $\U$set with its scope -- see \S\ref{S:NAFco}) if and only if $\U b$ is a $\U$member of the family of standard sets $\U$containing $\U a$, and $\U\nu(\U a)\in\nu(b)$ if and only if $\U a$ is a $\U$member of the set of standard members of $\U b$.\pa

\item Let $X$ be a Banach (resp.\ Hilbert) space. Then on $\U X$ there is a $\U \bR$-valued norm (resp.\ inner product and norm). In a well-known construction, one takes the space of elements of $\U X$ with {\em limited} norm and quotient it by the subspace of infinitesimal vectors (i.e.\ those with infinitesimal norm), (and denote equivalence classes modulo the infinitesimals by $[\,\,]$), to get an ordinary Banach (resp.\ Hilbert) space, called the {\em \NS\ hull} of $X$. Denote it by $\cH(X)$. It is complete, provided we have $\aleph_0^+$-saturation.%
\footnote{Indeed, Suppose $[\U x_n]$ is a Cauchy sequence in $\cH(X)$. For any $m\in\bN$, choose a ball $B_m=\{\U x\in\U X\,|\,\|\U x-\U a\|<1/m\}$ of radius $1/m$ which contains all $\U x_n$ from some $n$ on. The $B_m$ are internal sets and any finite number of them intersect (since they all contains all $\U x_n$ from some $n$ onward). By $\aleph_0^+$-Saturation they all intersect, and if $\U x_0$ is a member of the intersection, then clearly $[\U x_n]\to[\U x_0]$.}
The $\nu:H\to\U H$, followed by the factoring by the infinitesimals, defines a (norm-preserving) embedding of $X$ into $\cH(X)$.\pa

(In general, this embedding is not onto -- $\cH(X)$ has much more elements that $X$. For example, if $X=\ell^p$ then $\cH(X)$ consists of equivalence classes of ``sequences'' defined on $\U\bN$, and a sequence may have infinitesimal entries but finite nonzero norm, or can be supported on unlimited $\U$numbers, and then it is not in the image of $X$.)\pa

As for $\U\U X$, we have here, similarly to what we saw in $\bR$, the
$\U$infinitesimal vectors -- vectors with $\U$infinitesimal norm, which are part of the $\U\U$vactors with ``weakly infinitesimal'' norm. Call the latter ``weakly infinitesimal'' $\U\U$vactors. Factoring the space of $\U\U$vactors whose norm is bounded by a ($\nu\circ\nu$ of a) standard number by the ``weakly infinitesimal'' ones one again gets an ordinary Banach space (resp.\ Hilbert space) $\cH^{(2)}(X)$.\pa

Now, $\nu,\U\nu:\U X\to\U\U X$ will define two embeddings of $\cH(X)$ into $\cH^{(2)}(X)$, which by abuse of language we also refer to as $\nu$ and $\U\nu$.\pa

To see an example of what Prop.\ \ref{Prop:nunustar} would say here, take $X$ to be a Hilbert space. Suppose $[\U a]\in\cH(X)$ and $[\U b]$ belong to the orthogonal complement of $X$ in $\cH(X)$. Then for any (standard) $\eps>0$ the set of standard vectors $b$ with $|<b,\U a>|<\eps$ is everything. So by Prop.\ \ref{Prop:nunustar} $|<\U\nu(\U b),\nu(\U a)>|<\eps$. Consequently:\pa

If $[\U a]\in\cH(X)$ and $[\U b]$ belong to the orthogonal complement of $X$ in $\cH(X)$ then $[\U\nu(\U b)]$ is always orthogonal to $[\nu(\U a)]$. (In particular, $[\U\nu(\U b)]$ is orthogonal to $[\nu(\U b)]$.)
\end{itemize}
\bigskip

Another application of Prop.\ \ref{Prop:nunustar},
\begin{Prop}
Assume the \NA\ is $\aleph_0^+$-saturated.

Let $A$ be linearly ordered. (i.e.\ every $a,b\in A$ satisfy $a\ge b$ or $b\ge a$). Then a necessary and sufficient condition for $A$ to be well-ordered is that\pa

For every $\U a\in\U A$,\,\,$\U\nu(\U a)\le\nu(\U a)$.
\end{Prop}
\begin{Prf}
Suppose $A$ is well-ordered and $\U a\in\U A$. Let $c$ be the least (standard) element of $A$ satisfying $c>\U a$ --
if there is no $c$ with $c>\U a$ denote $c=\infty$. Then $\U a$ $\U$belongs to the set of standard elements $<c$, which is the set of standard elements $\le\U a$, hence by Prop.\ \ref{Prop:nunustar} $\U\nu(\U a)\le\nu(\U a)$.\pa

Suppose now that $A$ is not well-ordered, thus contains a strictly descending sequence $(a_n)$. By $\aleph_0^+$-saturation, there is a $\U a\in\U\{a_n\}$ smaller than all the $a_n$. The set of (standard) members of
$\{a_n\}$ that are $>\U a$ is the whole $\{a_n\}$, to which $\U a$ certainly $\U$belongs. Hence, by
Prop.\ \ref{Prop:nunustar}, $\U\nu(\U a)>\nu(\U a)$ and the condition is not satisfied.
\end{Prf}

This gives an immediate (\NS) proof to the fact that if some (standard) operation $\circ$ in $A$ satisfies
$$a'\ge a, b'\ge b\Rightarrow a'\circ a\ge b'\circ b$$
then if $S,S'\subset A$ are well-ordered then $S\circ S'$ is also well-ordered.\pa

Similar things can be proved for {\em partial} orders.

\appendix
\section{A Proof that an Infinite Set $B$ has $2^{2^{|B|}}$ Ultrafilters}\label{A:U}
For completeness, we give a proof for the well known
\begin{Prop}
The cardinality of the set $\be B$ of ultrafilters in an infinite set $B$ is $2^{2^{|B|}}$.
\end{Prop}

\begin{Prf}
Firstly, $|\be B|\le 2^{2^{|B|}}$ since $\be B\subset\cP\cP B$.

To prove the inverse inequality, it is enough to find a compact Hausdorff topological space $X$ of cardinality $2^{2^{|B|}}$ with a dense subset of cardinality $|B|$. Indeed, every point of $X$ is a limit of an ultrafilter in the dense subset, and an ultrafilter has a unique limit, hence there must be at least $2^{2^{|B|}}$ such ultrafilters.

We proceed to find such an $X$.

Denote by $\mathbf2=\{0,1\}$ the field, thus additive group, with two elements, endowed with the discrete topology.

Consider the group $\mathbf2^B$. By viewing the members of $\mathbf2^B$ as characteristic functions of subsets of $B$, we can identify $\mathbf2^B=\cP B$. $\mathbf2^B$ is a linear vector space over the field $\mathbf2$. As such, it has a (Hamel) basis $W$, and since the vector space is the set of sums of finite subsets of $W$, we have $|W|=|\mathbf2^B|=2^{|B|}$.

Let $H$ be the group of all homomorphisms $h:\mathbf2^B\to\mathbf2$. Since such a homomorphism is uniquely determined by giving its values on the members of $W$, which may be given arbitrarily, we may identify $H=\mathbf2^W$. The product topology here, which is also the weakest topology such that all evaluation maps at elements of $\mathbf2^B$ are continuous (hence does not depend on the choice of basis $W$), makes $H$ a compact (Hausdorff) topological group of cardinality $2^{2^{|B|}}$.

$B$ can be viewed as a subset of $H$, identifying every $b\in B$ with the homomorphism
$h:\mathbf2^B\to\mathbf2$ of evaluation at $b$. Let $<B>$ be the subgroup (also $\mathbf2$-linear subspace) of $H$ generated by this $B$. Since $<B>$ is just the set of
sums of finite subsets of $B$, we have $|<B>|=|B|$.

Thus, we shall be finished if we prove $<B>$ is dense in $H$.

Let $\BAR{<B>}$ be the closure of $<B>$ in $H$, and suppose $\BAR{<B>}\ne H$. Then there is a point $h_0\in H\setminus\BAR{<B>}$ and a neighborhood of $h_0$ disjoint from $\BAR{<B>}$. Thus there are $w_1,\ldots w_n\in W$ so that no element of $\BAR{<B>}$, in particular no element $\sum_\iota b_\iota\in<B>$, $(b_\iota)_\iota$ a finite subset of $B$, satisfies
$\sum_\iota w_1(b_\iota)=h_0(w_1),\ldots,\sum_\iota w_n(b_\iota)=h_0(w_n)$. This means that the set of all $\sum_\iota(w_1(b_\iota),\ldots,w_n(b_\iota))\in\mathbf2^n$ for all finite subsets $(b_\iota)_\iota$, is not the whole of $\mathbf2^n$. Since this set is a linear subspace of the $n$-dimensional vector space $\mathbf2^n$ over $\mathbf2$, and it is not the whole space, there is a non-zero linear functional which vanishes on it. Thus, there are $\eps_1,\ldots,\eps_n\in\mathbf2$, not all of them $0$, so that $\sum\eps_jw_j(b)=0$ for all $b\in B$, that is $\sum\eps_jw_j=0$, contradicting the linear independence of the basis $W$.
\end{Prf}

\end{document}